\documentclass[11pt]{amsart}
\usepackage[margin=1.15in]{geometry}
\usepackage[T1]{fontenc}
\usepackage{lmodern}
\usepackage[expansion=false]{microtype}
\usepackage{tikz}
\usepackage{amsmath,amssymb,amsthm,mathtools}
\usepackage{xcolor}
\usepackage[colorlinks=true,linkcolor=blue!50!black,citecolor=blue!50!black,urlcolor=blue!50!black]{hyperref}

\newtheorem{theorem}{Theorem}

\newtheorem*{theorem*}{Theorem}
\newtheorem{thm}{Theorem}[section]
\newtheorem{lemma}[thm]{Lemma}
\newtheorem{prop}[thm]{Proposition}

\theoremstyle{definition}
\newtheorem{defn}[thm]{Definition}
\theoremstyle{remark}
\newtheorem{remark}[thm]{Remark}

\newcommand{\N}{\mathfrak{N}}
\newcommand{\C}{\mathbb{C}}
\newcommand{\R}{\mathbb{R}}
\newcommand{\Z}{\mathbb{Z}}
\newcommand{\CH}{\mathbf{H}^2_{\C}}
\newcommand{\Ginf}{\Gamma_{\infty}}
\newcommand{\ip}[2]{\langle #1,#2\rangle}
\newcommand{\tst}{\tau_{*}}
\newcommand{\tcrit}{t_{c}}
\renewcommand{\Im}{\operatorname{Im}}
\renewcommand{\Re}{\operatorname{Re}}
\newcommand{\cl}{\overline}

\newcommand{\hatx}{\widehat{x}}

\hypersetup{pdftitle={A direct geometric proof of the Goldman--Parker theorem},pdfauthor={Krishnendu Gongopadhyay and Anna Pratoussevitch},pdfsubject={A direct geometric proof using a wall and a trapping region},pdfkeywords={Complex hyperbolic geometry, ideal triangle groups, Cygan metric, trapping regions, Goldman--Parker theorem}}
\title[The Goldman--Parker Conjecture]{A Short Proof of the Goldman--Parker Conjecture}
\author[K.~Gongopadhyay]{Krishnendu Gongopadhyay}
\address{Indian Institute of Science Education and Research (IISER) Mohali, Knowledge City,
S.A.S.\ Nagar, Sector 81, P.O. Manauli, 140306, India}
\email{krishnendu@iisermohali.ac.in}

\author[A.~Pratoussevitch]{Anna Pratoussevitch}
\address{Department of Mathematical Sciences, University of Liverpool, Peach Street,
Liverpool L69~7ZL, UK}
\email{annap@liverpool.ac.uk}

\date{November 25, 2026}
\begin{document}

\begin{abstract}
Goldman and Parker conjectured that the three order-two complex side
reflections of an ideal triangle with Cartan angular invariant $A$
define a discrete and faithful representation of $\Z/2 \ast  \Z/2 \ast  \Z/2$
in $PU(2,1)$ whenever $\tan^2 A\le 125/3$. Schwartz first proved
the conjecture using rigorous computer-assisted estimates and later
gave an analytic proof based on three loxodromic $\R$-spheres.

We give a substantially shorter proof of this conjecture using a simpler geometric
construction. For parameters just below the critical value, we use
the Cygan bisector of the two fixed points of the triple product of
the reflections. This wall and one of its cusp translates bound a
strip from which we construct a trapping region. At the critical
value $\tan^2 A=125/3$, the two fixed points merge, and the same
construction applies to the limiting wall.
\end{abstract}
\maketitle

\section{Introduction}

Let $\Gamma=\Z/2 \ast \Z/2 \ast \Z/2$, with standard generators $\varepsilon_0,\varepsilon_1,\varepsilon_2$. An ideal triangle in the complex hyperbolic plane, with distinct vertices and Cartan angular invariant $|A|<\pi/2$, determines a representation
\[
\rho_A\colon\Gamma\longrightarrow PU(2,1)
\]
by sending $\varepsilon_j$ to the order-two complex reflection $I_j$ in the complex geodesic through the other two vertices. Up to conjugacy and complex conjugation, the representation is determined by the absolute value of the Cartan angular invariant $A$. We use the parameter
\[
t=\tan|A|,\qquad 0\le t<\infty.
\]
Goldman and Parker \cite{GP} proved discreteness and faithfulness for $t^2\le35$ and showed that $t^2\le125/3$ is necessary. Schwartz proved that this necessary condition is sufficient \cite{S-Ann,S-GT}. The statement is often called the Goldman--Parker conjecture, after its formulation in~\cite{GP}; since Schwartz's proof we refer to it as the Goldman--Parker theorem. We give a direct geometric proof of the sufficiency statement, including the critical parameter.

\begin{theorem}[Goldman--Parker; Schwartz \cite{S-Ann,S-GT}]\label{thm:GP}
The representation $\rho_A$ is a discrete embedding for every $\tan^2A\le125/3$.
\end{theorem}

\subsection{Idea of the proof}
The proof treats the two types of elements of the cusp subgroup
separately. Write $\Gamma_\infty=\langle I_0,I_1\rangle$ and
$T=I_0I_1$. Every element of $\Gamma_\infty$ is either $T^n$ or
$T^nI_0$. The Ford ball is already disjoint from its images under
the reflected elements $T^nI_0$. For the translations, we use the
Cygan bisector of the two fixed points of $g=TI_2$, which we call
the \emph{wall}. This wall and its $T^{-1}$-translate bound an
open strip. An exact sum-of-squares identity shows that the strip
is disjoint from all its nontrivial $T$-translates.

Further identities for $g$, together with the positivity of a
two-by-two Hermitian matrix, show that $I_2$ sends the two open
regions outside the strip into the strip. Their images form two
caps. Adding these caps to the part of the Ford ball inside the
strip gives a region $U$. We call this construction of $U$ from
the wall, its translate and the Ford ball the \emph{wall
construction}. The region $U$ contains a neighbourhood of the
ball centre and satisfies $U\cap I_2(U)=\emptyset$ and
$I_2h(U)\subset U$ for every nontrivial $h\in\Gamma_\infty$.
These conditions imply discreteness and faithfulness. The region
$U$ may extend outside the Ford ball; a bound on the imaginary
part of the horizontal coordinate separates the caps from the
reflected balls.

At the critical parameter, the two fixed points become a single
parabolic fixed point. The formula defining the wall still has
a well-defined limit, and the same proof works with this point
on the boundary. The construction thus reduces the upper-range
proof to explicit identities for one wall and an estimate
separating the caps from the reflected balls. It avoids the
need to prove embeddedness and determine the intersections
of the three loxodromic $\mathbb R$-spheres used in \cite{S-GT}.

\subsection{The cusp subgroup and the initial ball}
Choose Heisenberg coordinates with the vertex fixed by $I_0$ and $I_1$ at infinity. Then
\[
\Ginf=\langle I_0,I_1\rangle\cong\Z/2\ast\Z/2
\]
is infinite dihedral. Put $T=I_0I_1$. In the normalization of Section~\ref{sec:setup}, the reflection $I_2$ exchanges infinity with
\[
o=(c,0),\qquad c=\tfrac12(1+it),
\]
and preserves the Ford isometric sphere $S_t$ with centre $o$ and radius $r=\tfrac12\sqrt{1+t^2}$. Write $B_t$ for its open ball. For the proof of Theorem~\ref{thm:GP}, it is enough to separate this ball from its nontrivial cusp-group translates for $t^2\le36$. Proposition~\ref{prop:initial36} does so by the triangle inequality and a short sum-of-squares identity. We also record the exact threshold, as a separate result. Put
\begin{equation}\label{eq:tau}
\tau_*=63-12w\approx37.478768,
\qquad w^3+3w=16,\quad w>0,
\qquad t_*:=\sqrt{\tau_*},\quad \tcrit:=\sqrt{125/3}.
\end{equation}

\begin{theorem}\label{thm:B}
The representation $\rho_A$ is a discrete embedding whenever $t^2\le\tau_*$. Moreover,
\[
B_t\cap T(B_t)=\emptyset
\quad\Longleftrightarrow\quad t^2\le\tau_*.
\]
For $t<t_*$ the corresponding closed balls are disjoint; at $t=t_*$ they meet at exactly one point.
\end{theorem}

The exact threshold describes the limit of the original ball
construction: beyond it, neighbouring ball translates overlap.
Its proof is in Appendix~\ref{app:sharp-ball} and is independent
of the proof of the full Goldman--Parker theorem. The wall
construction described above overcomes this limitation by using
a trapping region that may extend outside the Ford ball.
Theorem~\ref{thm:wall} verifies this construction for
$36\le t^2\le125/3$, including the critical parameter.
Together with Proposition~\ref{prop:initial36}, it proves
Theorem~\ref{thm:GP}.

\subsection{Relation to earlier work}\label{sec:related-work}
The Goldman--Parker family is a basic example of the deformation problem for discrete subgroups of $PU(2,1)$. At $A=0$, the representation preserves a totally real hyperbolic plane. Requiring the generators to remain complex reflections and their pairwise products to remain parabolic gives the one-parameter family considered here. Goldman's monograph \cite{G} supplies the geometric background, and Schwartz's survey \cite{S-ICM} places the family among complex hyperbolic triangle-group deformations.

Schwartz's first proof~\cite{S-Ann} uses rigorous computer-assisted
estimates. Numerical experiments guided his second proof~\cite{S-GT},
which establishes the required estimates analytically. He constructs
three loxodromic $\R$-spheres, one for each generating reflection,
and proves that each sphere is embedded and that each pair meets in
two arcs arranged in a T.

Our construction also uses the fixed points
of a triple product, but the argument simplifies by treating the two
types of cusp elements separately. For $t<\tcrit$, the Cygan bisector
of the fixed points of $g=TI_2$ gives a wall through the cusp. This
wall and its $T^{-1}$-translate bound a strip whose nontrivial
$T$-translates are disjoint. For the reflected elements $T^nI_0$,
a height estimate separates the caps from the reflected Ford balls.
At $t=\tcrit$, the two fixed points merge, and the limiting Hermitian
form defines a wall for which the same argument applies.
See~\cite{S-Deg} for the geometry at the critical parameter.

The two proofs also differ in the range of parameters they cover. The
loxodromic $\R$-sphere argument of~\cite{S-GT} is carried out for
$35\le t^2<125/3$. Below this range it relies on the theorem of
Goldman and Parker~\cite{GP}, and for the critical parameter it refers
to~\cite{S-Deg} and to Falbel and Parker~\cite{FP}. Schwartz notes that
computer pictures suggest his construction works for smaller parameters
as well, but the required technical estimates are established only on
this interval. The argument given here is self-contained over the whole
range $0\le t^2\le125/3$: the Ford ball covers $t^2\le36$, and the wall
construction covers $36\le t^2\le125/3$. At the critical parameter the
wall is defined by the limiting Hermitian form, and the identities of
Lemma~\ref{lem:four} hold there without change, so this case needs no
separate argument.

Table~\ref{tab:proof-comparison} summarizes the geometric
constructions and the main ideas in these proofs.

\begin{table}[htbp]
\centering
\small
\setlength{\tabcolsep}{5pt}
\renewcommand{\arraystretch}{1.25}
\begin{tabular}{@{}p{0.20\textwidth}
                  p{0.32\textwidth}
                  p{0.40\textwidth}@{}}
\hline
\textbf{Proof}
& \textbf{Geometric construction}
& \textbf{Main Calculation}
\\
\hline

Goldman--Parker \cite{GP}
&
Dirichlet regions for two dihedral subgroups preserving
complex geodesics.
&
Separation of the excluded regions, reduced to calculations
in hyperbolic discs. This gives the sufficient range
$t^2\le35$.
\\[0.7em]

Schwartz \cite{S-Ann}
&
A dented Clifford torus obtained by replacing two triangular
pieces with six hybrid sectors.
&
Placement and separation estimates for the hybrid sectors
and their images, established by rigorous computer-assisted
calculations using interval arithmetic.
\\[0.7em]

Schwartz \cite{S-GT}
&
Three loxodromic $\R$-spheres built from cones ruled by
arcs of $\R$-circles.
&
Embeddedness of the spheres and control of their pairwise
intersections, established through analytic estimates.
\\[0.7em]

Present paper
&
A wall and its translate bounding a strip, together with
the part of the Ford ball inside the strip and two caps.
&
Exact Hermitian-form identities give ordered wall translates
and cap inclusions. A bound on $\Im z$ separates the caps
from the reflected balls.
\\
\hline
\end{tabular}
\caption{Geometric constructions and verification steps
in the Goldman--Parker family.}
\label{tab:proof-comparison}
\end{table}

Goldman and Parker \cite{GP-D} study the cusp geometry of dihedral groups generated by inversions in asymptotic complex hyperplanes and construct Dirichlet regions in horospherical coordinates. Their setting includes $\Ginf$. Parker \cite{P-Ford} studies Ford isometric spheres and restrictions on their radii in discrete groups. We use the explicit Heisenberg action of $\Ginf$ and the Cygan inversion formula to compare one Ford ball with all its cusp-group translates.

The separation and trapping criteria are applications of ping-pong arguments; see Maskit \cite{Maskit} and Serre \cite{Se}. The reduced-word argument gives faithfulness, while an open set disjoint from each of its nontrivial images gives discreteness. Kalane and Parker \cite{KP} also use Heisenberg geometry to construct disjoint regions and give sufficient conditions for two unipotent parabolic maps to generate a free discrete subgroup of $PU(2,1)$.

The even subgroup of $\Gamma$ is generated by
\[
a=I_0I_1,\qquad b=I_1I_2,\qquad ab=I_0I_2.
\]
These three elements are unipotent parabolic. Thus this family lies in the two-generator parameter space of Parker and Will \cite{PW}, who construct a complex hyperbolic Riley slice using Ford domains and relate a boundary group to a spherical CR uniformisation of the Whitehead link complement. We retain the three reflections and use the fixed points of their triple product to define the wall.

Falbel and Koseleff \cite{FK} construct a four-dimensional family of discrete ideal-triangle-group representations using antiholomorphic generating involutions. Their generators have a different type from the complex reflections considered here. Pratoussevitch \cite{Prat} develops trace formulas for words in general complex hyperbolic triangle groups, giving algebraic obstructions arising from elliptic words. The present construction supplies the geometric separation needed for sufficiency in the Goldman--Parker family.

\subsection*{Organisation}
Section~\ref{sec:setup} fixes the coordinates and collects the generator formulas. Section~\ref{sec:criterion} gives the separation criteria, and Section~\ref{sec:proofA} gives the short initial ball estimate. Section~\ref{sec:trapping} gives the trapping criterion. Section~\ref{sec:wall} records the critical symmetry relation, constructs the wall and proves Theorem~\ref{thm:GP}. Appendix~\ref{app:sharp-ball} proves the exact Ford-ball threshold.

\section{Heisenberg coordinates and the normal form}\label{sec:setup}
The material in this section is mostly standard. We keep it for fixing conventions and colleacting useful formulas. For more information see \cite{G}, \cite{Pn}.

\subsection{Heisenberg space and the Cygan metric}
Let $\N=\C\times\R$ with the group law
\[
(z_1,t_1) \circ (z_2,t_2)=\bigl(z_1+z_2,\;t_1+t_2+2\Im(z_1\bar z_2)\bigr),\qquad (z,t)^{-1}=(-z,-t).
\]
The \emph{Cygan metric} is $\rho(p,q)=\|q^{-1} \circ p\|$, where $\|(z,t)\|=(|z|^4+t^2)^{1/4}$; explicitly,
\begin{equation}\label{eq:cygan}
\rho\bigl((z_1,t_1),(z_2,t_2)\bigr)^4=|z_1-z_2|^4+\bigl(t_1-t_2+2\Im(z_1\bar z_2)\bigr)^2 .
\end{equation}
This is a metric \cite{KR}. It is invariant under left translations $L_q(p)=q \circ p$ and under the rotations $(z,t)\mapsto(e^{i\theta}z,t)$, and the antiholomorphic map $(z,t)\mapsto(\bar z,-t)$ is also an isometry. We write $\pi(z,t)=z$ for vertical projection. Then $S^3=\N\cup\{\infty\}$.

\subsection{The Siegel model}
On $\C^3$ consider the Hermitian form
\[
\ip{u}{v}=u_1\bar v_3+u_2\bar v_2+u_3\bar v_1,
\]
which has signature $(2,1)$ and is linear in the first variable.
The complex hyperbolic plane $\CH$ is realised as the set of negative
complex lines in $\C^3$, and its ideal boundary as the set of null
complex lines. Every negative line has a representative
$(w_1,w_2,1)$, so this realisation identifies $\CH$ with the domain
\[
\bigl\{(w_1,w_2)\in\C^2:
       2\Re w_1+|w_2|^2<0\bigr\}.
\]
This unbounded domain is the \emph{Siegel model}.

Its ideal boundary is identified with $\N\cup\{\infty\}$.
The \emph{standard lifts} of $p=(z,t)\in\N$ and of $\infty$ are
\[
\widehat p=\bigl(-|z|^2+it,\ \sqrt2\,z,\ 1\bigr),
\qquad \widehat\infty=(1,0,0).
\]
These are null vectors. For $p=(z_1,t_1)$ and $q=(z_2,t_2)$,
a direct computation gives
\begin{equation}\label{eq:hermcygan}
\begin{gathered}
\ip{\widehat p}{\widehat q}
=-|z_1-z_2|^2
+i\bigl(t_1-t_2+2\Im(z_1\bar z_2)\bigr),\\
|\ip{\widehat p}{\widehat q}|=\rho(p,q)^2,
\qquad
\ip{\widehat p}{\widehat\infty}=1.
\end{gathered}
\end{equation}
The third coordinate of a vector $v$ is
$\ip{v}{\widehat\infty}$.

A positive vector $n$, meaning $\ip nn>0$, determines the complex
geodesic consisting of the negative lines in $n^\perp$.
The complex reflection with polar vector $n$ is represented by
\[
I_n(u)=-u+2\,\frac{\ip un}{\ip nn}\,n.
\]
This matrix lies in $SU(2,1)$ and is self-adjoint:
\[
\ip{I_nu}{v}=\ip{u}{I_nv}.
\]
Its projective action fixes the complex geodesic pointwise,
together with its boundary chain
$\{x\in\partial\CH:\ip{\widehat x}{n}=0\}$.

\subsection{The normal form}

\begin{lemma}\label{lem:normal}
Every ideal triangle with $|A|<\pi/2$ is equivalent, under $PU(2,1)$ and complex conjugation, to the triangle with vertices
\[
p_1=(0,0),\qquad p_0=(1,t),\qquad p_2=\infty,\qquad t=\tan|A|\ge0 .
\]
\end{lemma}

\begin{proof}
$PU(2,1)$ is doubly transitive on $S^3$, so we may assume $p_2=\infty$ and $p_1=(0,0)$. The stabiliser of these two points contains the dilations $(z,t)\mapsto(\lambda z,\lambda^2t)$ and the rotations $(z,t)\mapsto(e^{i\theta}z,t)$. Write $p_0=(z_0,t_0)$. If $z_0=0$, the three points lie on the chain $\{0\}\times\R\cup\{\infty\}$, so $|A|=\pi/2$. Otherwise we may normalise $z_0=1$.

By \eqref{eq:hermcygan}, $\ip{\widehat p_0}{\widehat p_1}=-1+it_0$ and $\ip{\widehat p_1}{\widehat\infty}=\ip{\widehat\infty}{\widehat p_0}=1$. Hence
\[
A=\arg\bigl(-\ip{\widehat p_0}{\widehat p_1}\ip{\widehat p_1}{\widehat p_2}\ip{\widehat p_2}{\widehat p_0}\bigr)=-\arctan t_0 .
\]
The isometry $(z,t)\mapsto(\bar z,-t)$ fixes $p_1$ and $p_2$ and changes the sign of $t_0$. It conjugates complex reflections to complex reflections, and conjugation preserves discreteness and faithfulness.
\end{proof}

We include the subscript $t$ when the dependence on the parameter needs to be stated. From now on we fix $t\ge0$ and the vertices of Lemma~\ref{lem:normal}. Let $I_j$ be the complex reflection fixing $p_{j-1}$ and $p_{j+1}$ (indices mod $3$).

\begin{lemma}\label{lem:gens}
Put $c=\tfrac12(1+it)$, $o=(c,0)\in\N$ and $r=|c|=\tfrac12\sqrt{1+t^2}$.
\begin{enumerate}
\item $I_0(z,s)=(-z,s)$ and $I_1(z,s)=(2-z,\ s+4\Im z)$. Their polar vectors are $n_0=(0,1,0)$ and $n_1=(-\sqrt2,1,0)$.
\item $I_2=I_{n_2}$ with $n_2=(0,1+it,\sqrt2)$. Its fixed chain $C_2$ passes through $p_0$ and $p_1$; it is the ellipse $\{(z,s):|z-c|=r,\ s=-2\Im(\bar c z)\}$ with centre $o$.
\item $T:=I_0I_1$ is the left translation $L_{(-2,0)}$, i.e.\ $T(z,s)=(z-2,\ s+4\Im z)$, and $T^n(z,s)=(z-2n,\ s+4n\Im z)$.
\item $\Ginf:=\langle I_0,I_1\rangle=\{T^n,\ T^nI_0:n\in\Z\}$ is infinite dihedral. All its elements are Cygan isometries fixing $\infty$, and they act on $\pi$ by $z\mapsto z-2n$ and $z\mapsto -z-2n$ respectively.
\end{enumerate}
\end{lemma}

\begin{proof}
(1) For $n=(0,1,0)$ we have $\ip{\widehat p}{n}=\sqrt2 z$, and $I_n(\widehat p)=(|z|^2-is,\sqrt2z,-1)$ is a multiple of the standard lift of $(-z,s)$. The formula for $I_1$ is the conjugate of $I_0$ by $L_{(1,0)}$. It is checked in the same way from the polar vector $(-\sqrt2,1,0)$, which is orthogonal to $\widehat\infty$ and to $\widehat{(1,s)}$ for every $s$.

(2) We have $\ip{\widehat p_1}{n_2}=0$, and $\ip{\widehat p_0}{n_2}=\sqrt2(-1+it)+\sqrt2(1-it)=0$. The description of $C_2$ follows from \eqref{eq:hermcygan}: the chain is the left translate by $o$ of the horizontal circle $|z|=r$.

(3) Compose the two formulas in (1), and compare with $(-2,0)*(z,s)=(z-2,s+4\Im z)$; iterate.

(4) $T$ has infinite order, so the dihedral group $\langle I_0,I_1\rangle$ is infinite. $I_0$ is a rotation and $T$ is a left translation, so both are Cygan isometries.

\end{proof}

\begin{lemma}[See Lemma 4.6 in \cite{Pn}]\label{lem:inversion}
$I_2(\infty)=o$, and $\rho(I_2x,o)\,\rho(x,o)=r^2$ for every $x\in\N\setminus\{o\}$. More generally,
\begin{equation}\label{eq:inversion2}
\rho(I_2x,I_2y)=\frac{r^2\,\rho(x,y)}{\rho(x,o)\,\rho(y,o)}\qquad(x,y\neq o).
\end{equation}
Consequently $I_2$ preserves the Cygan sphere $S=\{x:\rho(x,o)=r\}\supset C_2$ and exchanges the open ball $B=\{x:\rho(x,o)<r\}$ with the exterior $\{x:\rho(x,o)>r\}\cup\{\infty\}$.
\end{lemma}

\begin{proof}
Since $\ip{\widehat\infty}{n_2}=\sqrt2$ and $\ip{n_2}{n_2}=1+t^2=4r^2$,
\[
I_2\widehat\infty=\Bigl(-1,\ \tfrac{2\sqrt2(1+it)}{1+t^2},\ \tfrac{4}{1+t^2}\Bigr)=r^{-2}\,\widehat o .
\]
As $I_2$ is an involution, $I_2\widehat o=r^2\widehat\infty$. For $x\in\N$ write $I_2\widehat x=\lambda_x\,\widehat{I_2x}$. Comparing third coordinates,
\[
\lambda_x=\ip{I_2\widehat x}{\widehat\infty}=\ip{\widehat x}{I_2\widehat\infty}=r^{-2}\ip{\widehat x}{\widehat o},
\]
which is nonzero for $x\ne o$, and $|\lambda_x|=\rho(x,o)^2/r^2$. By \eqref{eq:hermcygan},
\[
\rho(I_2x,I_2y)^2=\frac{|\ip{I_2\widehat x}{I_2\widehat y}|}{|\lambda_x||\lambda_y|}=\frac{|\ip{\widehat x}{\widehat y}|}{|\lambda_x||\lambda_y|}=\frac{r^4\rho(x,y)^2}{\rho(x,o)^2\rho(y,o)^2},
\]
which is \eqref{eq:inversion2}. Taking $y\to\infty$ (or using $I_2\widehat o=r^2\widehat\infty$ directly) gives $\rho(I_2x,o)\rho(x,o)=r^2$.
\end{proof}

Translate the centre $o$ to the origin by $L_{o^{-1}}$. If $(\zeta,v)$ are the resulting coordinates, the fixed chain of $I_2$ is $|\zeta|=r$, $v=0$, and reflection in this chain is
\begin{equation}\label{eq:centredreflection}
(\zeta,v)\longmapsto
\left(\frac{r^2\zeta}{|\zeta|^2-iv},
      -\frac{r^4v}{|\zeta|^4+v^2}\right).
\end{equation}
This follows by using the polar vector $(r^2,0,1)$ in the reflection formula of \S\ref{sec:setup}.

$S$ is the Ford isometric sphere of $I_2$ relative to $\infty=p_2$. It is a spinal sphere (the boundary of a bisector) with poles $(c,\pm r^2)$, and $C_2$ is its equator.

\begin{lemma}[Centres]\label{lem:centres}
For $n\in\Z$,
\[
T^n o=(c-2n,\ 2nt),\qquad T^nI_0\,o=(-c-2n,\ -2nt),
\]
\[
\rho(o,T^no)^2=4|n|\sqrt{n^2+t^2},\qquad \rho(o,T^nI_0o)^2=(2n+1)^2+t^2 .
\]
\end{lemma}

\begin{proof}
$T^n=L_{(-2n,0)}$, so $T^no=(c-2n,-4n\Im\bar c)=(c-2n,2nt)$. Similarly $I_0o=(-c,0)$ and $T^nI_0o=(-c-2n,-2nt)$. Next,
\[
(T^no)^{-1}*o=(2n,\ -2nt+2\Im(c(\bar c-2n)))=(2n,-4nt),
\qquad
(T^nI_0o)^{-1}*o=(2(c+n),0).
\]
Taking norms gives the two distances.
\end{proof}

\begin{lemma}[The antiholomorphic symmetry]\label{lem:J}
Let $J(z,s)=(1-\bar z,\ t-s-2\Im z)$. Then $J$ is an antiholomorphic Cygan isometry and an involution. It fixes $\infty$ and $o$, it exchanges $p_0$ and $p_1$, and it satisfies $JI_2J=I_2$, $JI_0J=I_1$, hence $JTJ=T^{-1}$ and $J(T^nI_0)J=T^{-n-1}I_0$.
\end{lemma}

\begin{proof}
$J=L_{(1,t)}\circ R$, where $R(z,s)=(-\bar z,-s)$ is an antiholomorphic Cygan isometry. A direct computation gives $J^2=\mathrm{id}$, $J(0,0)=(1,t)$ and $J(c,0)=(c,0)$. $J$ maps the chain $C_2$ to a chain through $p_0$ and $p_1$; since a chain is determined by two of its points, $J(C_2)=C_2$. Hence $JI_2J$ is a holomorphic involution fixing $C_2$ pointwise, i.e.\ $I_2$. Similarly $J$ maps $\{0\}\times\R$ onto $\{1\}\times\R$, so $JI_0J=I_1$. The last two identities follow.
\end{proof}

\begin{remark}[The dihedral subgroup]
In these coordinates $I_0$, $I_1$ and hence $\Ginf$ are independent of $t$. The identities $I_0TI_0=T^{-1}$ and $I_1TI_1=T^{-1}$ express the reversibility of the translation $T$. Every element of $\Ginf$ has a unique expression as $T^n$ or $T^nI_0$. If $K\subset\N$ is compact, only finitely many $h\in\Ginf$ satisfy $hK\cap K\ne\emptyset$: the real parts of horizontal coordinates in $K$ are bounded, whereas their images under $T^n$ and $T^nI_0$ have real parts $\Re z-2n$ and $-\Re z-2n$, respectively. Thus $\Ginf$ acts properly discontinuously on $\N$; finite stabilizers are allowed.
\end{remark}

\section{The one-sphere criterion}\label{sec:criterion}

The criteria in this section are elementary formulations of the ping-pong argument for free products. The reduced-word argument proves faithfulness. Discreteness follows from the existence of a nonempty open set disjoint from all its nontrivial translates.

All complementary components in the following statement are taken in $\partial\CH=\N\cup\{\infty\}$.

\begin{thm}[One-sphere criterion]\label{thm:criterion}
Suppose a smooth embedded sphere $S\subset\N$ bounds a relatively compact region $U\subset\N$ and an exterior region $V$ containing infinity. Assume
\[
I_2(S)=S,\qquad I_2(U)=V,
\qquad U\cap h(U)=\emptyset\quad(h\in\Ginf\setminus\{1\}).
\]
Then the natural homomorphism
\[
\Ginf*\langle I_2\rangle\longrightarrow PU(2,1)
\]
is injective and has discrete image. In particular $\rho_A$ is a discrete embedding.
\end{thm}
\begin{proof}
The open set $h(U)$ is disjoint from $U$, so it is also disjoint from $S=\partial U$: every open neighbourhood of a point of $S$ meets $U$. Thus $h(U)\subset V$ for every nontrivial $h\in\Ginf$. Put
\[
X_2=U,\qquad X_1=\bigcup_{h\in\Ginf\setminus\{1\}}h(U).
\]
These sets are disjoint, $h(X_2)\subset X_1$ for $h\ne1$, and $I_2(X_1)\subset X_2$.

To prove discreteness, we choose an open set disjoint from all its nontrivial images. Choose $R_0$ with $\pi(\cl U)\subset D(c,R_0)$, where $c=(1+it)/2$. The set
\[
\Omega=\{(z,s):|\Im z|>t/2+R_0\}
\]
is $\Ginf$-invariant and disjoint from every $h(\cl U)$. For $Y>t/2+R_0+1$, take
\[
W=\{(z,s):|z-(1/4+iY)|<1/8\}\subset\Omega.
\]
The horizontal actions in Lemma~\ref{lem:gens} show that $h(W)\cap W=\emptyset$ for every nontrivial $h\in\Ginf$.

A reduced word involving $I_2$ first takes $W$ into $V$, then by $I_2$ into $U$, and subsequently alternates between $X_2$ and $X_1$. Its final image lies in $\bigcup_hh(U)$, disjoint from $W$. Words in $\Ginf\setminus\{1\}$ also send $W$ to a set disjoint from itself. Thus every nontrivial reduced word sends $W$ to a set disjoint from it, which proves injectivity.

If nonidentity elements of the image tended to the identity in $PU(2,1)$, their images of any fixed point of $W$ would eventually lie in $W$. This contradicts the preceding conclusion. Hence the image is discrete.
\end{proof}

\begin{thm}[Functional separation criterion]\label{thm:functional}
Let $P\colon\N\to\R$ satisfy
\begin{enumerate}
\item $P(o)<0$, and $P(x)$ and $P(I_2x)$ have opposite signs or are both zero whenever $x\ne o$;
\item the set $\{P\le0\}$ is bounded;
\item $P(hx)>0$ whenever $P(x)<0$ and $h\in\Ginf\setminus\{1\}$.
\end{enumerate}
Then $\rho_A$ is a discrete embedding. No regularity or topological assumption on the zero set is required.
\end{thm}
\begin{proof}
Put $U=\{P<0\}$, $V=\{P>0\}$ and $X_1=\bigcup_{h\ne1}h(U)$, $X_2=U$. The hypotheses give $X_1\subset V$, $h(X_2)\subset X_1$ for $h\ne1$, and $I_2(X_1)\subset X_2$. Boundedness allows us to choose the horizontal region $\Omega$ and open set $W$ as in Theorem~\ref{thm:criterion}, with the latter disjoint from every $h\{P\le0\}$. Every $\Ginf$-translate of $W$ lies in $V$. The reduced-word argument in that theorem therefore applies without change, and proves injectivity and discreteness.
\end{proof}

\subsection{Trapping regions}\label{sec:trapping}

The following combination argument applies to a region that need not be exchanged with its full complement by $I_2$. We use it for the construction in Section~\ref{sec:wall}.

\begin{thm}[Trapping criterion]\label{thm:D}
Let $U\subset\N$ satisfy
\begin{enumerate}
\item[\textup{(T1)}] $U\cap I_2(U)=\emptyset$;
\item[\textup{(T2)}] $U$ contains an open neighbourhood of $o$;
\item[\textup{(T3)}] $I_2h(U)\subset U$ for every $h\in\Ginf\setminus\{1\}$.
\end{enumerate}
Then $\Ginf*\langle I_2\rangle\to PU(2,1)$ is injective with discrete image. In particular $\rho_A$ is a discrete embedding.
\end{thm}

\begin{proof}
By \textup{(T1)}, the function $P\colon\N\to\R$ defined by $P=-1$ on $U$, $P=1$ on $I_2U\cap\N$ and $P=0$ elsewhere is well defined. We check the hypotheses of Theorem~\ref{thm:functional}.

(1) $P(o)=-1$. For $x\ne o$, $x\in U$ if and only if $I_2x\in I_2U$, and $x\in I_2U$ if and only if $I_2x\in U$. Hence $P(x)$ and $P(I_2x)$ have opposite signs or both vanish.

(2) If $U\supset\{\rho(x,o)<\epsilon\}$, then Lemma~\ref{lem:inversion} gives $I_2U\supset\{\rho(x,o)>r^2/\epsilon\}$. Hence $\{P\le0\}=\N\setminus I_2U$ is bounded.

(3) If $P(x)<0$, then $x\in U$ and, by \textup{(T3)}, $I_2(hx)\in U$. Thus $hx\in I_2U$ and $P(hx)=1>0$.
\end{proof}

Condition \textup{(T3)} says that every translate $hU$, $h\ne1$, lies in $I_2U$. We do not require $U$ and $I_2U$ to fill the boundary: the sign function is zero elsewhere. The region $U$ may extend outside $B$. If $U\subset B$, condition \textup{(T1)} follows from Lemma~\ref{lem:inversion}.

\section{The Ford ball construction}\label{sec:proofA}

We first prove discreteness and faithfulness for $t^2\le36$
using the Ford ball. The triangle inequality separates the ball
from its reflected images and its distant translation images.
Only the neighbouring translates require a further estimate,
which follows from a short sum-of-squares identity. This gives
the initial range needed before applying the wall construction.

We apply Theorem~\ref{thm:criterion} with $S=S_t$ and $U=B_t$. Invariance and exchange of components follow from Lemma~\ref{lem:inversion}. Each translate $h(B_t)$ is the Cygan ball of radius $r$ centred at $h(o)$. Recall that $4r^2=1+t^2$.

\begin{lemma}[Far translates]\label{lem:far}
If $t^2\le65$ and $h\in\Ginf\setminus\{1,T,T^{-1}\}$, then $h(B)\cap\cl B=\emptyset$.
\end{lemma}

\begin{proof}
By the triangle inequality it suffices to show $\rho(o,ho)^2\ge4r^2=1+t^2$; we use Lemma~\ref{lem:centres}.
\begin{itemize}
\item For $h=T^nI_0$: $\rho(o,ho)^2=(2n+1)^2+t^2\ge1+t^2$.
\item For $h=T^n$ with $|n|\ge2$: $\rho(o,ho)^2\ge 8\sqrt{4+t^2}$. The inequality $8\sqrt{4+t^2}\ge1+t^2$ is equivalent to $t^4-62t^2-255\le0$. For $\tau=t^2\in[0,65]$, one has $\tau^2-62\tau-255\le3\tau-255\le-60$, as required. \qedhere
\end{itemize}
\end{proof}

\begin{lemma}[Intersection of neighbouring balls]\label{lem:key}
Let $t\ge0$. Then $B\cap T(B)=\emptyset$ if and only if $t^2\le\tst=63-12w$, where $w^3+3w=16$. If $t^2\le\tst$, then moreover $T(B)\cap\cl B=\emptyset$ and $B\cap T(\cl B)=\emptyset$. When $t^2=\tst$, the closed balls meet in exactly one point; when $t^2<\tst$ they are disjoint.
\end{lemma}

\begin{proof}
The sharp calculation is given in Appendix~\ref{app:sharp-ball}. The proof of the full Goldman--Parker theorem uses only the coarser bound in Proposition~\ref{prop:initial36} below.
\end{proof}

\begin{prop}[The initial range]\label{prop:initial36}
If $t^2\le36$, then
\[
h(B)\cap\cl B=\emptyset \qquad(h\in\Ginf\setminus\{1\}).
\]
Consequently $\rho_A$ is a discrete embedding in this range.
\end{prop}
\begin{proof}
Lemma~\ref{lem:far} handles all translates except $T$ and $T^{-1}$. If $t^2\le15$, then
\[
\rho(o,To)^2=4\sqrt{1+t^2}\ge1+t^2=4r^2,
\]
so the triangle inequality gives the required separation.

For $15\le t^2\le36$, apply the Cygan isometry
\[
\Lambda=L_{(1,-2t)}\circ L_{o^{-1}},\qquad
\Lambda(o)=(1,-2t),\quad\Lambda(To)=(-1,2t).
\]
At a point $x=(X+iy,v)$, let $D_+(x)$ and $D_-(x)$ be the fourth powers of its distances from these two centres. Formula~\eqref{eq:cygan} gives
\begin{equation}\label{eq:ball-average}
\frac{D_+(x)+D_-(x)}2
=(1+y^2)^2+4(t+y)^2+X^4+(6+2y^2)X^2+v^2.
\end{equation}
An intersection of the two balls, with at least one ball open, would therefore imply $Q_t(y)<0$, where
\begin{equation}\label{eq:Qt}
Q_t(y):=16\bigl[(1+y^2)^2+4(t+y)^2\bigr]-(1+t^2)^2.
\end{equation}
The identity
\begin{equation}\label{eq:initial-SOS}
\begin{split}
Q_t(y)={}&16(y^2-3)^2+192\left(y+\frac t3\right)^2\\
&+39+(36-t^2)\left(t^2-\frac{14}{3}\right)
\end{split}
\end{equation}
shows that $Q_t(y)\ge39$ for every $y\in\R$ whenever $15\le t^2\le36$, since the last product is then nonnegative. (The identity in fact gives $Q_t>0$ for $14/3\le t^2\le36$.) This proves both $T(B)\cap\cl B=\emptyset$ and $B\cap T(\cl B)=\emptyset$; applying $T^{-1}$ gives the corresponding assertion for $T^{-1}(B)$. Theorem~\ref{thm:criterion} now applies.
\end{proof}

\begin{proof}[Proof of Theorem~\ref{thm:B}]
Let $t^2\le\tst$. By Lemma~\ref{lem:key} the open ball $T(B)$ is disjoint from $\cl B$, and $T^{-1}(B)\cap\cl B=T^{-1}\bigl(B\cap T(\cl B)\bigr)=\emptyset$. Since $\tst<65$, Lemma~\ref{lem:far} handles all other $h\in\Ginf\setminus\{1\}$, so the disjointness hypothesis of Theorem~\ref{thm:criterion} holds. Theorem~\ref{thm:criterion} applies.
\end{proof}

\begin{remark}\label{rem:triangle}
The triangle inequality alone gives only $t^2\le15$. Indeed, $\rho(o,To)^2=4\sqrt{1+t^2}\ge1+t^2$ if and only if $t^2\le15$. The sum-of-squares identity~\eqref{eq:initial-SOS} extends this to $t^2\le36$, which is sufficient for the full theorem. The exact threshold in Lemma~\ref{lem:key} is proved separately in Appendix~\ref{app:sharp-ball}. 
\end{remark}

\section{A wall through the cusp and the full threshold}\label{sec:wall}

Put
\[
g:=TI_2=I_0I_1I_2,\qquad g':=JgJ=T^{-1}I_2 .
\]
Goldman and Parker showed that $g$ is loxodromic for $t^2<125/3$ and parabolic at $t^2=125/3$ \cite{GP}. Indeed, $\operatorname{tr}g=(17i-t)/(t-i)$, and Goldman's discriminant \cite[\S6.2]{G}
\[
f(\tau)=|\tau|^4-8\Re(\tau^3)+18|\tau|^2-27
\]
takes the value $f(\operatorname{tr}g)=1024(125-3t^2)/(1+t^2)^3$. In this section we prove the full threshold with a trapping region in the sense of Theorem~\ref{thm:D}. We do not need to calculate the overlap $B\cap T(\cl B)$. The isometric ball already gives disjoint regions for the reflections $T^nI_0$, so we keep it for those elements. For the translations $T^n$ it is replaced by a strip between a hypersurface $\Sigma$ through $\infty$ and its translate $T^{-1}\Sigma$. For $t<\tcrit$ the hypersurface $\Sigma$ is the Cygan bisector of the two fixed points of $g$. Three facts are needed:
\begin{itemize}
\item $T$ moves $\Sigma$ strictly to one side. The sum-of-squares identity~\eqref{eq:wallidentity} proves this.
\item $g$ moves $\Sigma$ to the same side. For $t<\tcrit$ the translates $g^n\Sigma$ are pairwise disjoint; at $t=\tcrit$ they all pass through the parabolic fixed point $q$ of $g$ and meet only there.
\item $I_2=T^{-1}g$ maps one side of $\Sigma$ off itself. We prove this using a positive Hermitian form that splits into blocks on the eigenspaces of $I_2$.
\end{itemize}
All three properties hold for every $0<t\le\tcrit$. 

\subsection{A rational parametrisation}

\begin{lemma}\label{lem:param}
The function $\mu\mapsto(\mu-2)(2\mu+1)^2/\bigl(3(\mu+2)\bigr)$ is an increasing bijection of $[2,7]$ onto $[0,125/3]$. Let $0<t\le\tcrit$, let $\mu\in(2,7]$ be defined by
\[
t^2=\frac{(\mu-2)(2\mu+1)^2}{3(\mu+2)},
\]
and put
\[
\nu:=\sqrt{\frac{\mu^2-4}{3}},\qquad \xi:=\sqrt{\frac{(7-\mu)(\mu+1)}{12}} .
\]
Then
\begin{gather}
t=\frac{2\mu+1}{\mu+2}\,\nu,\qquad t-\nu=\frac{\mu-1}{\mu+2}\,\nu,\qquad 2\nu-t=\frac{3}{\mu+2}\,\nu,\label{eq:param1}\\
1+\nu^2=\frac{\mu^2-1}{3},\qquad r^2=\frac{(\mu-1)^2(\mu+1)}{3(\mu+2)}=(1+\nu^2)\,\frac{t-\nu}{\nu},\label{eq:param2}\\
\xi^2=\frac{t}{\nu}-\frac{(t-\nu)^2}{4}.\label{eq:param3}
\end{gather}
In particular $0<t-\nu<\nu$ and $\nu\ge3t/5$. Moreover $\xi=0$ exactly when $t=\tcrit$, and then $\nu=\sqrt{15}$.
\end{lemma}

\begin{proof}
On $[2,7]$ the factors $(\mu-2)/(\mu+2)$ and $(2\mu+1)^2$ are nonnegative and increasing, and the values at $2$ and $7$ are $0$ and $5\cdot225/27=125/3$. Since $\nu^2=(\mu-2)(\mu+2)/3$, the definition of $\mu$ reads $t^2=(2\mu+1)^2\nu^2/(\mu+2)^2$, which gives \eqref{eq:param1}. Next,
\[
1+t^2=\frac{3(\mu+2)+(\mu-2)(2\mu+1)^2}{3(\mu+2)}=\frac{4(\mu-1)^2(\mu+1)}{3(\mu+2)},
\]
and $4r^2=1+t^2$. The second expression for $r^2$ follows from \eqref{eq:param1} and $1+\nu^2=(\mu^2-1)/3$. Finally
\[
\frac t\nu-\frac{(t-\nu)^2}{4}=\frac{2\mu+1}{\mu+2}-\frac{(\mu-1)^2(\mu-2)}{12(\mu+2)}
=\frac{-\mu^3+4\mu^2+19\mu+14}{12(\mu+2)}=\frac{(7-\mu)(\mu+1)}{12}.
\]
The inequalities follow from \eqref{eq:param1}, since $0<(\mu-1)/(\mu+2)<1$ and $(2\mu+1)/(\mu+2)=2-3/(\mu+2)\le5/3$.
\end{proof}

Throughout this section $\mu$, $\nu$ and $\xi$ are the quantities of Lemma~\ref{lem:param}, and
\[
\omega:=\xi+\tfrac i2(t-\nu),\qquad \varkappa:=\frac{t+\nu}{\nu}.
\]

\begin{lemma}[Centred coordinates]\label{lem:centred}
For $x=(z,s)\in\N$ put
\[
\zeta(x):=z-c,\qquad \varpi(x):=-\ip{\hatx}{\widehat o}=|z-c|^2-i\bigl(s+2\Im(z\bar c)\bigr).
\]
The pair $(\zeta,\varpi)$ determines $x$. If $x\ne o$, then
\begin{equation}\label{eq:I2centred}
\zeta(I_2x)=\frac{r^2\,\zeta(x)}{\varpi(x)},\qquad \varpi(I_2x)=\frac{r^4}{\varpi(x)} ,
\end{equation}
and for every $x\in\N$
\begin{equation}\label{eq:Tcentred}
\zeta(T^{-1}x)=\zeta(x)+2,\qquad \varpi(T^{-1}x)=\varpi(x)+4\zeta(x)+8c .
\end{equation}
Consequently a point $x\in\N\setminus\{o\}$ is fixed by $g$ if and only if
\begin{equation}\label{eq:fixedeq}
r^2\zeta=\varpi\,(\zeta+2),\qquad r^4=\varpi\,(\varpi+4\zeta+8c),
\end{equation}
and in that case $g\hatx=-r^{-2}\varpi(x)\,\hatx$.
\end{lemma}

\begin{proof}
The formula for $\varpi$ is \eqref{eq:hermcygan}. In the proof of Lemma~\ref{lem:inversion} we saw that $I_2\hatx=\lambda_x\widehat{I_2x}$ with $\lambda_x=r^{-2}\ip{\hatx}{\widehat o}=-\varpi(x)/r^2$, and that $I_2\widehat o=r^2\widehat\infty$. Hence
\[
\varpi(I_2x)=-\lambda_x^{-1}\ip{I_2\hatx}{\widehat o}=-\lambda_x^{-1}\ip{\hatx}{I_2\widehat o}=-\lambda_x^{-1}r^2=\frac{r^4}{\varpi(x)} .
\]
The third coordinate of $I_2\hatx=-\hatx+\tfrac{1}{2r^2}\ip{\hatx}{n_2}n_2$ is $\lambda_x$, so $\ip{\hatx}{n_2}=\sqrt2\,r^2(1+\lambda_x)$. Its second coordinate is therefore $-\sqrt2z+\sqrt2c(1+\lambda_x)$, and
\[
\zeta(I_2x)=\frac{-z+c(1+\lambda_x)}{\lambda_x}-c=-\frac{\zeta(x)}{\lambda_x}=\frac{r^2\zeta(x)}{\varpi(x)} .
\]
This is \eqref{eq:I2centred}; it is the reflection \eqref{eq:centredreflection}. For \eqref{eq:Tcentred}, $T$ preserves the form and maps standard lifts to standard lifts, and $To=(c-2,2t)$. Thus, by \eqref{eq:hermcygan},
\[
\varpi(T^{-1}x)=-\ip{\hatx}{\widehat{To}}=|\zeta+2|^2-i\bigl(s-2t+2\Im(z(\bar c-2))\bigr).
\]
Since $\Im z=\Im\zeta+t/2$, the right-hand side equals $\varpi+4\Re\zeta+4+4i\Im\zeta+4it=\varpi+4\zeta+8c$. Finally $gx=x$ means $I_2x=T^{-1}x$, which is \eqref{eq:fixedeq} by \eqref{eq:I2centred} and \eqref{eq:Tcentred}. In that case $I_2\hatx=\lambda_x\widehat{I_2x}$, and $T$ maps the standard lift of $I_2x$ to the standard lift of $TI_2x=x$. Hence $g\hatx=\lambda_x\hatx=-r^{-2}\varpi(x)\,\hatx$.
\end{proof}

\subsection{Fixed points and the critical orbit}
The following explicit relation identifies the boundary points retained by the wall construction at the critical parameter.

\begin{prop}\label{prop:criticalobstruction}
Let $t=\tcrit=5\sqrt{15}/3$, and let $I_2$, $J$ and $T$ denote the maps of Lemmas~\ref{lem:gens} and~\ref{lem:J} at this parameter. Put
\[
q=\left(-\tfrac12+\tfrac{\sqrt{15}}2i,\ \tfrac{4\sqrt{15}}3\right).
\]
Then
\begin{equation}\label{eq:criticalorbit}
I_2q=Jq=T^{-1}q
 =\left(\tfrac32+\tfrac{\sqrt{15}}2i,-\tfrac{2\sqrt{15}}3\right).
\end{equation}
Moreover $q\in S\cap T(S)$.
\end{prop}
\begin{proof}
The formulas for $J$ and $T^{-1}$ give the displayed coordinates directly. For $I_2$, translate $o$ to the origin as in \eqref{eq:centredreflection}. The point $q$ then has coordinates
\[
\left(-1-\tfrac{\sqrt{15}}3i,\ \tfrac{8\sqrt{15}}3\right),\qquad r^2=\tfrac{32}{3}.
\]
Equation~\eqref{eq:centredreflection}, followed by translation back, gives the same image in \eqref{eq:criticalorbit}. For the last assertion, the centred coordinates of $q$ give $\rho(q,o)^4=\bigl(1+\tfrac{15}9\bigr)^2+\tfrac{64\cdot15}{9}=\tfrac{1024}9=r^4$, so $q\in S$. Then $T^{-1}q=I_2q\in S$, that is, $q\in T(S)$.
\end{proof}

The point $q$ is fixed by $g=TI_2$, which is parabolic at the critical parameter \cite{GP}; see Proposition~\ref{prop:fixedpoints}. The wall construction below retains $q$ and $Jq$ on the boundary of its trapping region.

\begin{prop}[The fixed points of $g$]\label{prop:fixedpoints}
Let $0<t\le\tcrit$. Then
\begin{equation}\label{eq:omega}
|\omega|^2=\frac t\nu,\qquad |\omega\pm1|^2=\varkappa\pm2\xi,\qquad r^2=(1-i\nu)(\omega-1)(\bar\omega+1).
\end{equation}
The points
\begin{equation}\label{eq:qpm}
q_\pm:=\Bigl(-\tfrac12\mp\xi+\tfrac{i\nu}2,\ \ \tfrac{t+\nu}2\pm(2\nu-t)\,\xi\Bigr)
\end{equation}
are fixed by $g$. Their standard lifts satisfy $g\widehat q_\pm=\lambda_\pm\widehat q_\pm$ with
\[
|\lambda_+|=\frac{|\omega+1|}{|\omega-1|},\qquad |\lambda_-|=\frac{|\omega-1|}{|\omega+1|} .
\]
If $t<\tcrit$, then $q_+\ne q_-$, and $q_+$ is the attracting and $q_-$ the repelling fixed point of $g$. If $t=\tcrit$, then $q_+=q_-$ is the parabolic fixed point $q$ of $g$, the point of Proposition~\ref{prop:criticalobstruction}.
\end{prop}

\begin{proof}
The first identity in \eqref{eq:omega} follows from \eqref{eq:param3}. The second follows because $\Re\omega=\xi$. For the third,
\[
(\omega-1)(\bar\omega+1)=|\omega|^2-1+(\omega-\bar\omega)=\frac{t-\nu}{\nu}+i(t-\nu)=\frac{t-\nu}{\nu}\,(1+i\nu),
\]
so multiplying by $1-i\nu$ and using \eqref{eq:param2} gives the result.

We compute the centred coordinates of $q_\pm$. Since $c=\tfrac12(1+it)$,
\[
\zeta(q_+)=-(1+\omega),\qquad \zeta(q_-)=\bar\omega-1 ,
\]
so $\Re\varpi(q_\pm)=|\omega\pm1|^2$. For the imaginary part, $2\Im(z\bar c)=\Im z-t\,\Re z$, which at $q_\pm$ equals $\tfrac12(t+\nu)\pm t\xi$. Adding the vertical coordinate in \eqref{eq:qpm} gives $t+\nu\pm2\nu\xi=\nu(\varkappa\pm2\xi)$. Hence, by \eqref{eq:omega},
\begin{equation}\label{eq:varpiq}
\varpi(q_\pm)=(1-i\nu)\,|\omega\pm1|^2 .
\end{equation}

Now we verify \eqref{eq:fixedeq} for $q_+$. The first equation reads $-r^2(1+\omega)=(1-i\nu)(1+\omega)(1+\bar\omega)(1-\omega)$, which is the third identity in \eqref{eq:omega}. Given the first equation, $r^4/\varpi=r^2(\zeta+2)/\zeta=r^2(\omega-1)/(\omega+1)$, and by the complex conjugate of \eqref{eq:omega} this is $(1+i\nu)|\omega-1|^2$. The second equation in \eqref{eq:fixedeq} therefore reads
\[
(1-i\nu)|\omega+1|^2-(1+i\nu)|\omega-1|^2=4(1+\omega)-4(1+it)=4(\omega-it).
\]
By \eqref{eq:omega} the left-hand side is $4\xi-2i\nu\varkappa=4\xi-2i(t+\nu)$, and the right-hand side is $4\xi+2i(t-\nu)-4it$. They agree. For $q_-$ the first equation is again the third identity in \eqref{eq:omega}, and the second is the complex conjugate of the identity just displayed.

By Lemma~\ref{lem:centred} and \eqref{eq:varpiq}, $|\lambda_\pm|=\sqrt{1+\nu^2}\,|\omega\pm1|^2/r^2$, and $r^2=\sqrt{1+\nu^2}\,|\omega-1|\,|\omega+1|$ by \eqref{eq:omega}. This gives the stated eigenvalue moduli. If $t<\tcrit$, then $\xi>0$, so the two points are distinct and $|\lambda_+|>1>|\lambda_-|$. If $t=\tcrit$, then $\xi=0$, $\nu=\sqrt{15}$ and $t=5\sqrt{15}/3$, and \eqref{eq:qpm} is the point $q$.
\end{proof}

\begin{remark}
The common height $\Im z=\nu/2$ of the two fixed points provides a parameter for the family: by \eqref{eq:param1}, $t/\nu=2-3/(\mu+2)$ increases from $5/4$ to $5/3$. In particular $\nu<t<2\nu$. Both inequalities are used in Lemma~\ref{lem:four} below.
\end{remark}

\subsection{The Wall}

Recall that $T=I_0I_1$, and put $g=TI_2=I_0I_1I_2$.
The two fixed points of $g$ determine the wall.

The two fixed points of $g=TI_2$ determine the wall. We first express
the construction using their standard lifts, then use coordinate
formulas that remain meaningful when the fixed points come together
at the critical parameter.

\begin{defn}\label{def:wall}
Let $0<t\le\tcrit$. For $t<\tcrit$, let $q_+$ and $q_-$ be the
attracting and repelling fixed points of $g$, and put
\[
\mathbf e:=\widehat q_++\widehat q_-,
\qquad
\mathbf d:=\frac{\widehat q_--\widehat q_+}{2\xi}.
\]
Thus $\mathbf e$ is the sum of the fixed-point lifts, while
$\mathbf d$ is their scaled difference. Their coordinate
expressions are
\[
\begin{aligned}
\mathbf d
&=\bigl(1-i(2\nu-t),\ \sqrt2,\ 0\bigr),\\
\mathbf e
&=\bigl(-\tfrac12(1+\nu^2)-2\xi^2+i(t+\nu),\
        \sqrt2\,(-1+i\nu),\ 2\bigr).
\end{aligned}
\]
We use these expressions also at $t=\tcrit$, where they define
the vectors by continuity.

For $u\in\C^3$, define
\[
\mathcal F(u)
:=\Re\Bigl(\ip u{\mathbf d}\,
                 \overline{\ip u{\mathbf e}}\Bigr),
\qquad
\mathcal G(u)
:=\tfrac12|\ip u{\mathbf e}|^2
  +2\xi^2|\ip u{\mathbf d}|^2.
\]
Both are Hermitian forms, and $\mathcal G\ge0$. Their values
are multiplied by $|\lambda|^2$ when $u$ is replaced by
$\lambda u$, so the sign of $\mathcal F$ is well defined
on projective points. Put
\[
\mathfrak H_+:=\{\mathcal F>0\},\qquad
\mathfrak H_-:=\{\mathcal F<0\},\qquad
\Sigma:=\{\mathcal F=0\}
\]
as subsets of $\partial\CH$. The sets $\mathfrak H_\pm$ are open.
Since the third coordinate of $\mathbf d$ is zero,
$\ip{\widehat\infty}{\mathbf d}=0$, and hence $\infty\in\Sigma$.
\end{defn}

The identities
\[
\widehat q_\pm=\tfrac12\mathbf e\mp\xi\mathbf d
\]
give
\begin{equation}\label{eq:Fbisector}
\begin{aligned}
2\xi\,\mathcal F(u)
&=|\ip u{\widehat q_-}|^2-|\ip u{\widehat q_+}|^2,\\
\mathcal G(u)
&=|\ip u{\widehat q_+}|^2+|\ip u{\widehat q_-}|^2.
\end{aligned}
\end{equation}
Consequently, for $t<\tcrit$ and $x\in\N$,
\[
\mathcal F(\widehat x)
=\frac{\rho(x,q_-)^4-\rho(x,q_+)^4}{2\xi}.
\]
Thus $\Sigma\cap\N$ is the Cygan bisector of the two fixed
points, and $\mathfrak H_+\cap\N$ consists of the points closer
to $q_+$ than to $q_-$. The form $\mathcal G$ is the corresponding
sum of the fourth powers of the two distances.

At $t=\tcrit$, the fixed points coincide at $q$, and
$\mathbf e=2\widehat q$. The wall is then defined by
\[
\mathcal F(u)
=2\Re\Bigl(\ip u{\mathbf d}\,
                 \overline{\ip u{\widehat q}}\Bigr).
\]
The scaled difference $\mathbf d$ retains the direction in which
the fixed-point lifts separate. This gives the limiting wall
at the critical parameter.

\begin{lemma}\label{lem:Fcoords}
For $x=(z,s)\in\N$ put $\chi:=\Re z+\tfrac12$, $y:=\Im z$, and
\[
\mathsf n:=\chi^2+\bigl(y-\tfrac\nu2\bigr)^2+\xi^2,\qquad
\varsigma:=s-\tfrac t2-\nu\chi-y,\qquad
\ell:=2y+2\nu-t .
\]
Then $\ip{\hatx}{\mathbf d}=2\chi+i\ell$ and $\ip{\hatx}{\mathbf e}=-2\mathsf n+2i\varsigma$. Consequently
\begin{equation}\label{eq:Fcoords}
\mathcal F(\hatx)=2\varsigma\ell-4\chi\,\mathsf n,\qquad
\mathcal G(\hatx)=2(\mathsf n^2+\varsigma^2)+2\xi^2(4\chi^2+\ell^2).
\end{equation}
In the coordinates $(\chi,y,s)$ the maps $T$ and $J$ are
\begin{equation}\label{eq:TJcoords}
T(\chi,y,s)=(\chi-2,\ y,\ s+4y),\qquad J(\chi,y,s)=(2-\chi,\ y,\ t-s-2y).
\end{equation}
\end{lemma}

\begin{proof}
By the definition of the form, $\ip{\hatx}{\mathbf d}=\sqrt2z\cdot\sqrt2+\overline{\mathbf d_1}=2z+1+i(2\nu-t)$. Similarly
\[
\ip{\hatx}{\mathbf e}=2(-|z|^2+is)-2z(1+i\nu)-\tfrac12(1+\nu^2)-2\xi^2-i(t+\nu).
\]
Its real part is $-2|z|^2-2\Re z+2\nu y-\tfrac12(1+\nu^2)-2\xi^2=-2\mathsf n$, and its imaginary part is $2s-2y-2\nu\Re z-t-\nu=2\varsigma$. Formula \eqref{eq:TJcoords} follows from Lemma~\ref{lem:gens}(3) and Lemma~\ref{lem:J}.
\end{proof}

\subsection{The linear algebra behind the wall identities}

Two elementary linear-algebra facts explain the identities \eqref{eq:gmonotone}  and 
\eqref{eq:reflection-blocks} used to prove
\eqref{eq:I2certificate} for $I_2$. The first uses the attracting
and repelling eigenvectors of $g$, and the second uses the
eigenspace decomposition of the involution $I_2$.

\begin{lemma}[Eigenvector scaling]\label{lem:eigenvector-scaling}
Let $L$ preserve a Hermitian form of signature $(n,1)$, and let $v_+,v_-$ be nonzero null eigenvectors with eigenvalues $\lambda_+,\lambda_-$ satisfying
\[
|\lambda_+|=R>1,\qquad |\lambda_-|=R^{-1}.
\]
Put $A_\pm(u)=|\ip u{v_\pm}|^2$, $\Phi=A_--A_+$ and $\Psi=A_-+A_+$. Then
\begin{equation}\label{eq:eigenvector-scaling}
\Phi(L^{\pm1}u)
=\frac{R^2+R^{-2}}2\,\Phi(u)
\ \pm\ \frac{R^2-R^{-2}}2\,\Psi(u).
\end{equation}
Moreover, $\Psi(u)>0$ for every nonzero null vector $u$.
\end{lemma}
\begin{proof}
Unitarity gives $\ip{Lu}{v_\pm}=\ip u{L^{-1}v_\pm}=\bar\lambda_\pm^{-1}\ip u{v_\pm}$. Thus $A_+(Lu)=R^{-2}A_+(u)$ and $A_-(Lu)=R^2A_-(u)$, which proves the identity; applying $L^{-1}$ reverses the sign. The span of the two independent null eigenvectors has signature $(1,1)$, so its orthogonal complement is positive definite. A nonzero null vector cannot be orthogonal to both eigenvectors.
\end{proof}

\begin{lemma}[Averaging over an involution]\label{lem:involution-averaging}
Let $I$ be a linear involution of a complex vector space, and let $Q$ be a Hermitian form. Write $u=u_++u_-$ with $Iu_\pm=\pm u_\pm$. Then
\[
Q(u)+Q(Iu)=2Q(u_+)+2Q(u_-).
\]
If $I$ preserves a Hermitian form $H$, its two eigenspaces are $H$-orthogonal. Consequently $cH-Q-Q\circ I$ has no mixed terms between them.
\end{lemma}
\begin{proof}
Expand $Q(u_++u_-)+Q(u_+-u_-)$. The mixed terms cancel. Also $H(u_+,u_-)=H(Iu_+,Iu_-)=-H(u_+,u_-)$, proving orthogonality.
\end{proof}

\begin{lemma}[Four properties of the wall]\label{lem:four}
Let $0<t\le\tcrit$. In the following, $g$ and $I_2$ denote the matrices $TI_{n_2}$ and $I_{n_2}$.
\begin{enumerate}
\item[\textup{(a)}] For $x\in\N$, $\mathcal F(\widehat{Jx})=-\mathcal F(\widehat{Tx})$. Hence $J(\mathfrak H_\pm)=T^{-1}(\mathfrak H_\mp)$.
\item[\textup{(b)}] For $x\in\N$,
\begin{equation}\label{eq:wallidentity}
\mathcal F(\widehat{Tx})-\mathcal F(\hatx)=24(\chi-1)^2+24\Bigl(y+\frac{\nu}{2(\mu+2)}\Bigr)^2+\frac{8(\mu+1)^2}{\mu+2} .
\end{equation}
In particular $\mathcal F(\widehat{Tx})>\mathcal F(\hatx)$.
\item[\textup{(c)}] For every $u\in\C^3$,
\begin{equation}\label{eq:gmonotone}
(\varkappa^2-4\xi^2)\,\mathcal F(g^{\pm1}u)=(\varkappa^2+4\xi^2)\,\mathcal F(u)\pm2\varkappa\,\mathcal G(u),
\end{equation}
and $\varkappa^2-4\xi^2=|\omega^2-1|^2>0$. Consequently, if $t<\tcrit$, the sets $g^n\Sigma$, $n\in\Z$, are pairwise disjoint. If $t=\tcrit$, the parabolic fixed point $q$ of $g$ lies on every $g^n\Sigma$, and $g^m\Sigma\cap g^n\Sigma=\{q\}$ for $m\ne n$.
\item[\textup{(d)}] For every $u\in\C^3\setminus\{0\}$,
\begin{equation}\label{eq:I2certificate}
\mathcal F(u)+\mathcal F(I_2u)<\frac{4(1+\nu^2)}{r^2}\,\ip uu .
\end{equation}
In particular $\mathcal F(u)+\mathcal F(I_2u)<0$ for every null vector $u\ne0$.
\end{enumerate}
\end{lemma}

\begin{proof}
(a) and (b). Under $T$ the quantities of Lemma~\ref{lem:Fcoords} become $\mathsf n-4\chi+4$, $\varsigma+4y+2\nu$ and $\ell$, by \eqref{eq:TJcoords}. Under $J$ they become $\mathsf n-4\chi+4$, $-(\varsigma+4y+2\nu)$ and $\ell$, while $\chi$ becomes $2-\chi$. Substituting in \eqref{eq:Fcoords} gives
\[
\mathcal F(\widehat{Tx})=2(\varsigma+4y+2\nu)\ell-4(\chi-2)(\mathsf n-4\chi+4)=-\mathcal F(\widehat{Jx}).
\]
This proves (a); the statement about $\mathfrak H_\pm$ follows because $J$ is an involution. Subtracting \eqref{eq:Fcoords},
\begin{align*}
\mathcal F(\widehat{Tx})-\mathcal F(\hatx)
&=8\mathsf n+16\chi^2-48\chi+32+4(2y+\nu)(2y+2\nu-t)\\
&=24(\chi-1)^2+24y^2+8(2\nu-t)\,y+8+8\xi^2+2\nu^2+4\nu(2\nu-t).
\end{align*}
Completing the square in $y$ and using $2\nu-t=3\nu/(\mu+2)$ leaves the constant
\[
8+8\xi^2+2\nu^2+4\nu(2\nu-t)-\tfrac23(2\nu-t)^2
=8+(4\mu+2)+4(\mu-2)-\frac{2(\mu-2)}{\mu+2}=\frac{8(\mu+1)^2}{\mu+2}.
\]
Here $8\xi^2+2\nu^2=\tfrac23(6\mu+3)$, $4\nu(2\nu-t)=12\nu^2/(\mu+2)=4(\mu-2)$ and $\tfrac23(2\nu-t)^2=6\nu^2/(\mu+2)^2$.

\medskip (c) Suppose first $t<\tcrit$, and put $R=|\omega+1|/|\omega-1|$. By Proposition~\ref{prop:fixedpoints}, Lemma~\ref{lem:eigenvector-scaling} applies to $g$ with $v_\pm=\widehat q_\pm$. Equation~\eqref{eq:Fbisector} gives $\Phi=2\xi\mathcal F$ and $\Psi=\mathcal G$. Since
\[
R^{\pm2}=\frac{\varkappa\pm2\xi}{\varkappa\mp2\xi},\qquad
\frac{R^2+R^{-2}}2=\frac{\varkappa^2+4\xi^2}{\varkappa^2-4\xi^2},\qquad
\frac{R^2-R^{-2}}{4\xi}=\frac{2\varkappa}{\varkappa^2-4\xi^2},
\]
the scaling identity is exactly~\eqref{eq:gmonotone}. Its denominator is $|\omega-1|^2|\omega+1|^2>0$, because $\Im\omega>0$. The coefficients of both sides depend continuously on $\mu\in(2,7]$, so the identity also holds at $t=\tcrit$.

For a nonzero null vector $u$, Lemma~\ref{lem:eigenvector-scaling} gives $\mathcal G(u)>0$ before the critical parameter. At the critical parameter, $\xi=0$ and $\mathbf e=2\widehat q$, so
\[
\mathcal G(u)=2|\ip u{\widehat q}|^2.
\]
This vanishes only when $[u]=q$: in signature $(2,1)$, two orthogonal null vectors are proportional. Thus~\eqref{eq:gmonotone} sends $\mathcal F\ge0$ into $\mathcal F>0$, apart from the critical point $q$. Iteration shows that no positive iterate of a point of $\Sigma$ returns to $\Sigma$, except for $q$ at the critical parameter. This proves the asserted intersections. Finally $\mathcal F(\widehat q)=0$ and $gq=q$, so $q$ belongs to every $g^n\Sigma$ at that parameter.

\medskip \textup{(d)} Put $\mathbf n=(-\sqrt2\bar c,1,0)$. Then
$\ip{\mathbf n}{\mathbf n}=1$, and $\mathbf n$ is orthogonal
to $\widehat o$ and $\widehat\infty$. Moreover,
\[
I_2\mathbf n=-\mathbf n,\qquad
I_2\widehat o=r^2\widehat\infty,\qquad
I_2(r^2\widehat\infty)=\widehat o.
\]
Write
\[
u=\alpha\widehat o+\beta r^2\widehat\infty+\gamma\mathbf n,
\qquad \alpha_\pm=\alpha\pm\beta.
\]
Then $u=u_++u_-$, where
\[
u_+=\frac{\alpha_+}{2}
       (\widehat o+r^2\widehat\infty),\qquad
u_-=\frac{\alpha_-}{2}
       (\widehat o-r^2\widehat\infty)+\gamma\mathbf n,
\]
satisfy $I_2u_\pm=\pm u_\pm$. Since $I_2$ preserves the
ambient Hermitian form, these eigenspaces are orthogonal, and
\[
\ip uu
=\frac{r^2}{2}
  \bigl(|\alpha_+|^2-|\alpha_-|^2\bigr)+|\gamma|^2.
\]

Lemma~\ref{lem:involution-averaging} gives
\[
\mathcal F(u)+\mathcal F(I_2u)
=2\mathcal F(u_+)+2\mathcal F(u_-).
\]
Consequently,
\[
\frac{4(1+\nu^2)}{r^2}\ip uu
-\mathcal F(u)-\mathcal F(I_2u)
\]
has no mixed terms between the two eigenspaces. Its restriction
to the $+1$ eigenspace is a scalar form in $\alpha_+$, while
its restriction to the $-1$ eigenspace is a two-by-two
Hermitian form in $(\alpha_-,\gamma)$.

Set
\[
\vartheta=\varkappa+i(t+\nu),\qquad
\mathsf M=2r^2+1+\nu^2.
\]
To compute the two blocks, we use the pairings
\begin{align*}
\ip{\widehat o}{\mathbf d}&=2(1+i\nu),&
\ip{r^2\widehat\infty}{\mathbf d}&=0,&
\ip{\mathbf n}{\mathbf d}&=\sqrt2,\\
\ip{\widehat o}{\mathbf e}&=-2\vartheta,&
\ip{r^2\widehat\infty}{\mathbf e}&=2r^2,&
\ip{\mathbf n}{\mathbf e}&=-\sqrt2(2-i(t-\nu)).
\end{align*}
Together with the identities
\begin{equation}\label{eq:threesums}
\varkappa+\nu(t+\nu)=r^2+2(1+\nu^2),\qquad
r^2-\nu(t-\nu)=\frac{t-\nu}{\nu},
\end{equation}
which follow from~\eqref{eq:param2}, these give
\begin{align*}
2\mathcal F(u_+)
&=-4(1+\nu^2)|\alpha_+|^2,\\
2\mathcal F(u_-)
&=-4(r^2+1+\nu^2)|\alpha_-|^2-8|\gamma|^2
  -4\sqrt2\,\Re(\vartheta\alpha_-\bar\gamma).
\end{align*}
Combining these formulas with the expression for $\ip uu$,
we obtain
\begin{equation}\label{eq:reflection-blocks}
\begin{split}
\frac{4(1+\nu^2)}{r^2}\ip uu
-\mathcal F(u)-\mathcal F(I_2u)
={}&6(1+\nu^2)|\alpha_+|^2\\
&+\begin{pmatrix}\bar\alpha_-&\bar\gamma\end{pmatrix}
K\begin{pmatrix}\alpha_-\\\gamma\end{pmatrix},
\end{split}
\end{equation}
where
\[
K=\begin{pmatrix}
2\mathsf M&2\sqrt2\,\bar\vartheta\\
2\sqrt2\,\vartheta&4\mathsf M/r^2
\end{pmatrix}.
\]
The first diagonal entry of $K$ is positive, and
\[
\det K
=\frac8{r^2}
 \bigl(\mathsf M^2-r^2|\vartheta|^2\bigr).
\]
Put $a=t-\nu$. By~\eqref{eq:param2},
\[
r^2=(1+\nu^2)\frac a\nu,\qquad
\mathsf M=(1+\nu^2)\frac{2a+\nu}{\nu},\qquad
|\vartheta|^2
=(1+\nu^2)\frac{(a+2\nu)^2}{\nu^2}.
\]
Consequently,
\[
\begin{aligned}
\mathsf M^2-r^2|\vartheta|^2
&=\frac{(1+\nu^2)^2}{\nu^3}
  \bigl[\nu(2a+\nu)^2-a(a+2\nu)^2\bigr]\\
&=\frac{(1+\nu^2)^2}{\nu^3}(\nu^3-a^3)>0,
\end{aligned}
\]
because $\nu>0$ and $a<\nu$ by Lemma~\ref{lem:param}.
Thus $K$ is positive definite. Since
$(\alpha_+,\alpha_-,\gamma)\ne(0,0,0)$ whenever $u\ne0$,
the right-hand side of~\eqref{eq:reflection-blocks}
is strictly positive. This proves~\eqref{eq:I2certificate}.
\end{proof}

\begin{remark}\label{rem:wallmeaning}
Property \textup{(b)} says that the translates $T^n\Sigma$
are pairwise disjoint in $\N$ and occur in order: their
positive sides satisfy
\[
T^{n+1}(\mathfrak H_+)\subset T^n(\mathfrak H_+)
\qquad(n\in\Z).
\]
Property \textup{(c)} gives the corresponding inclusion
$g(\mathfrak H_+)\subset\mathfrak H_+$.
For $t<\tcrit$, the walls $g^n\Sigma$ are pairwise disjoint.
At $t=\tcrit$, they all contain the parabolic fixed point $q$,
and any two meet only at $q$.

The connection between \textup{(c)} and \textup{(d)} is
particularly useful. Property \textup{(d)} gives
$I_2(\mathfrak H_+)\subset\mathfrak H_-$.
Since $g=TI_2$, property \textup{(c)} also gives
$I_2(\mathfrak H_+)\subset T^{-1}(\mathfrak H_+)$.
Together these imply
\[
I_2(\mathfrak H_+)
\subset
\mathfrak H_-\cap T^{-1}(\mathfrak H_+),
\]
the open strip between $\Sigma$ and $T^{-1}\Sigma$
used below. Thus the two properties place the reflected
region on the required side of each wall.

At the critical parameter, the common point $q$ remains
on the boundary of the strip. The construction uses the
open sides, so this boundary contact does not affect
the argument.
\end{remark}

\subsection{The strip and the two caps}

Let
\[
\mathcal W:=\mathfrak H_-\cap T^{-1}(\mathfrak H_+),\qquad
\mathcal O:=\mathfrak H_+\cup T^{-1}(\mathfrak H_-),\qquad
\mathcal C:=I_2(\mathcal O),\qquad
\mathcal C^{*}:=I_2\bigl(\partial\CH\setminus\mathcal W\bigr).
\]
Thus $\mathcal W=\{x\in\N:\ \mathcal F(\hatx)<0<\mathcal F(\widehat{Tx})\}$ is the open strip between $\Sigma$ and $T^{-1}\Sigma$, the set $\mathcal O$ is the union of the two open sides of the strip, and $\mathcal C\subset\mathcal C^{*}$. Since $\infty\in\Sigma$, we have $\mathcal W\cup\mathcal O\subset\N$ and $o=I_2(\infty)\in\mathcal C^{*}$.

\begin{lemma}\label{lem:strip}
Let $0<t\le\tcrit$.
\begin{enumerate}
\item $\mathcal W\cap\mathcal O=\emptyset$, and $T^n(\mathcal W)\subset\mathcal O$ for every $n\ne0$.
\item $o\in\mathcal W$.
\item $\mathcal C\subset\mathcal W$.
\end{enumerate}
\end{lemma}

\begin{proof}
(1) The first statement follows from the definitions. Let $x\in\mathcal W$, and put $f_k:=\mathcal F(\widehat{T^kx})$ for $k\in\Z$. By Lemma~\ref{lem:four}(b) the sequence $f_k$ is strictly increasing, and $f_0<0<f_1$ by the definition of $\mathcal W$. If $n\ge1$, then
$f_n\ge f_1>0$, so $T^nx\in\mathfrak H_+$. If $n\le-1$, then $\mathcal F\bigl(\widehat{T(T^nx)}\bigr)=f_{n+1}\le f_0<0$, so $T^nx\in T^{-1}(\mathfrak H_-)$.

(2) At $o$ we found $\chi=1$, $\mathsf n=\varkappa$, $\varsigma=-(t+\nu)$ and $\ell=2\nu$. At $To=(c-2,2t)$ we have $\chi=-1$, $\mathsf n=\varkappa$, $\varsigma=t+\nu$ and $\ell=2\nu$. By \eqref{eq:Fcoords},
\[
\mathcal F(\widehat o)=-4\nu(t+\nu)-4\varkappa<0<4\nu(t+\nu)+4\varkappa=\mathcal F(\widehat{To}).
\]

(3) Let $x\in\mathfrak H_+$. Then $x\ne\infty$, and $x\ne o$ by (2), so $I_2x\in\N$. By Lemma~\ref{lem:four}(d), $\mathcal F(I_2\hatx)<-\mathcal F(\hatx)<0$, so $I_2x\in\mathfrak H_-$. Since $TI_2=g$, Lemma~\ref{lem:four}(c) gives
\[
(\varkappa^2-4\xi^2)\,\mathcal F(TI_2\hatx)=(\varkappa^2+4\xi^2)\,\mathcal F(\hatx)+2\varkappa\,\mathcal G(\hatx)>0,
\]
so $TI_2x\in\mathfrak H_+$. Hence $I_2(\mathfrak H_+)\subset\mathcal W$. By Lemma~\ref{lem:four}(a) and $JTJ=T^{-1}$,
\[
J(\mathcal W)=T^{-1}(\mathfrak H_+)\cap\mathfrak H_-=\mathcal W,\qquad T^{-1}(\mathfrak H_-)=J(\mathfrak H_+).
\]
As $J$ commutes with $I_2$, we get $I_2\bigl(T^{-1}(\mathfrak H_-)\bigr)=J\,I_2(\mathfrak H_+)\subset J(\mathcal W)=\mathcal W$.
\end{proof}

The set $\mathcal C$ consists of the two \emph{caps} $g^{-1}(\mathfrak H_-)=I_2T^{-1}(\mathfrak H_-)$ and $I_2(\mathfrak H_+)=J\,g^{-1}(\mathfrak H_-)$. For $t<\tcrit$ the first contains the repelling fixed point $q_-$ of $g$, the second contains the repelling fixed point $Jq_-$ of $g'=JgJ$, and their closures meet at $o$.

\begin{lemma}[The caps lie above the lower row of balls]\label{lem:caps}
Let $6\le t\le\tcrit$. Then $\mathcal C^{*}\subset\N$, and $\Im z>\tfrac1{12}$ at every point of $\mathcal C^{*}$.
\end{lemma}

\begin{proof}
By Lemma~\ref{lem:four}(a), $\partial\CH\setminus\mathcal W=\{\mathcal F\ge0\}\cup J\bigl(\{\mathcal F\ge0\}\cap\N\bigr)$. The maps $T^{-1}$ and $J$ preserve $\N$ and the function $\Im z$, $J$ commutes with $I_2$, and $I_2=T^{-1}g$. It therefore suffices to show that every point $x$ of $g\{\mathcal F\ge0\}$ lies in $\N$ and satisfies $\Im z>\tfrac1{12}$.

Let $x$ be such a point, so that $\mathcal F(g^{-1}\hatx)\ge0$. Then $x\ne\infty$, because $g^{-1}(\infty)=o$ and $\mathcal F(\widehat o)<0$. By Lemma~\ref{lem:four}(c), $(\varkappa^2+4\xi^2)\mathcal F(\hatx)\ge2\varkappa\,\mathcal G(\hatx)$. Put
\[
\mathsf a:=\frac{\varkappa^2+4\xi^2}{4\varkappa},\qquad \mathsf b:=\frac{\varkappa^2-4\xi^2}{4\varkappa},
\qquad\text{so that}\qquad \mathsf a+\mathsf b=\frac\varkappa2,\quad \mathsf a^2-\xi^2=\mathsf b^2 .
\]
By \eqref{eq:Fcoords} the inequality reads $\mathsf n^2+\varsigma^2+\xi^2(4\chi^2+\ell^2)\le\mathsf a\,(2\varsigma\ell-4\chi\mathsf n)$, that is,
\[
(\mathsf n+2\mathsf a\chi)^2+(\varsigma-\mathsf a\ell)^2\le\mathsf b^2\,(4\chi^2+\ell^2).
\]
Hence $\mathsf n+2\mathsf a\chi\le\mathsf b\,(2|\chi|+|\ell|)$. Since $2\mathsf b|\chi|-2\mathsf a\chi-\chi^2\le2(\mathsf a+\mathsf b)|\chi|-\chi^2\le(\mathsf a+\mathsf b)^2$ and $(\mathsf a+\mathsf b)^2-\xi^2=\tfrac14\varkappa^2-\xi^2=\mathsf b\varkappa$, the definition of $\mathsf n$ gives
\begin{equation}\label{eq:capbound}
\bigl(y-\tfrac\nu2\bigr)^2\le\mathsf b\,\bigl(\varkappa+|2y+2\nu-t|\bigr).
\end{equation}

For $6\le t\le\tcrit$, Lemma~\ref{lem:param} gives $\nu\ge18/5>7/2$. Also
\[
\varkappa=\frac{3(\mu+1)}{\mu+2}\le\frac83,\qquad
\mathsf b=\frac{(\mu-1)^3}{36(\mu+2)}\le\frac23,\qquad
0<2\nu-t=\sqrt{\frac{3(\mu-2)}{\mu+2}}<\frac43.
\]
The three expressions are increasing for $2<\mu\le7$; the formula for $\mathsf b$ follows by substituting $\varkappa$ and $\xi^2$ from Lemma~\ref{lem:param}.
Suppose $y\le1/12$ and put $h=1/12-y\ge0$. Then $\nu/2-y>5/3+h$ and $|2y+2\nu-t|\le3/2+2h$. Equation~\eqref{eq:capbound} would imply
\[
\left(\frac53+h\right)^2
<\frac23\left(\frac83+\frac32+2h\right)
=\frac{25}{9}+\frac43h,
\]
or $h^2+2h<0$, a contradiction. Hence $y>1/12$.
\end{proof}

\subsection{The trapping region and the proof of the full threshold}

Put
\[
U_\flat:=(B\cap\mathcal W)\cup\mathcal C .
\]

\begin{figure}[htbp]
\centering
\input{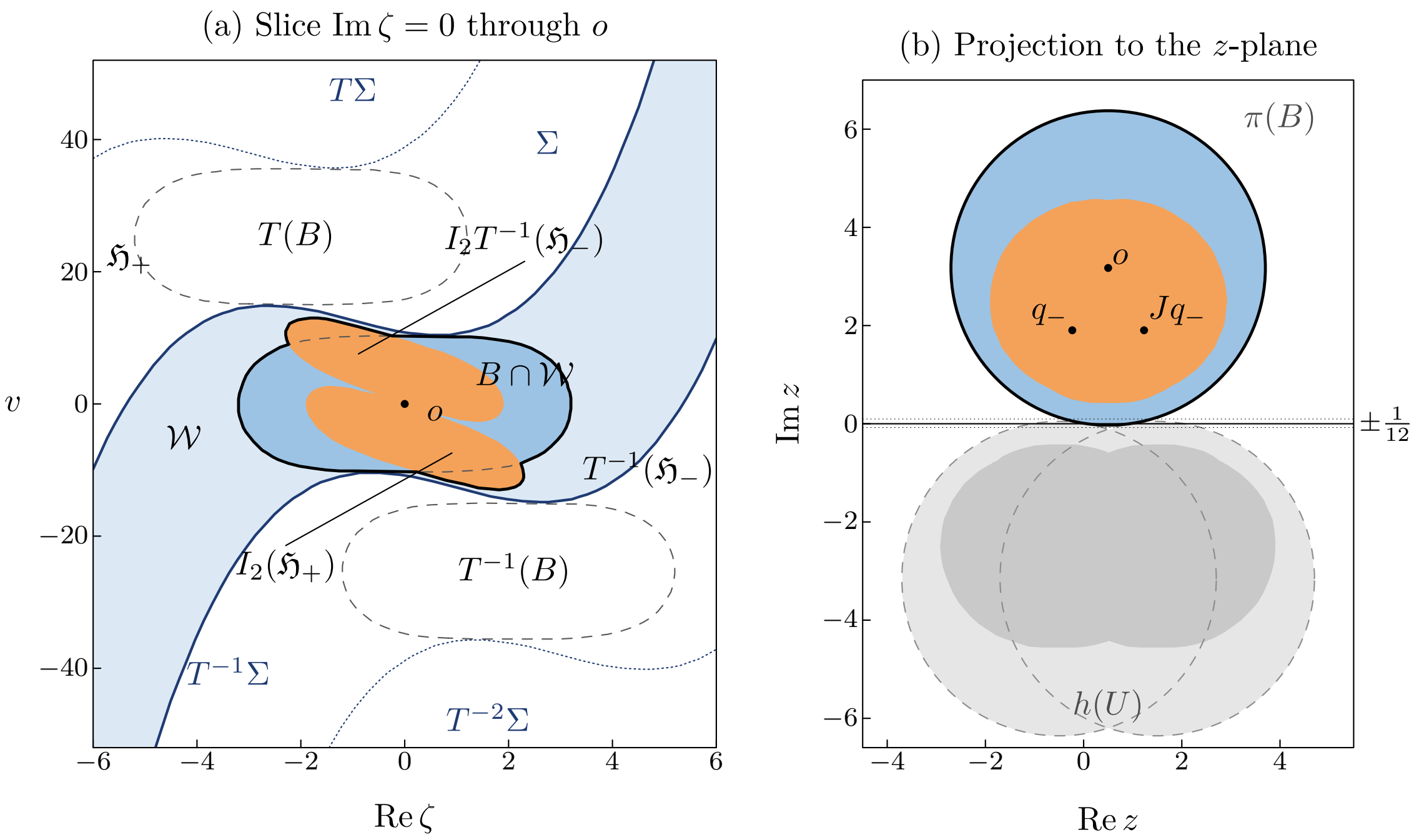}
\caption{The wall construction at $t^2=40$.
(a) The slice $\operatorname{Im}\zeta=0$ through $o$, in the centred
coordinates $\zeta=z-c$, $v=s+2\operatorname{Im}(z\bar c)$. The strip
$\mathcal W$ lies between $\Sigma$ and $T^{-1}\Sigma$, and the
translates $T^{\pm1}(B)$ lie outside it. The caps
$\mathcal C=I_2(\mathcal O)$ (orange) have closures meeting at $o$, and the
trapping region $U_\flat=(B\cap\mathcal W)\cup\mathcal C$ (outlined)
extends outside the Ford ball $B$ (dashed).
(b) Vertical projection to the $z$-plane. The caps lie in
$\operatorname{Im}z>1/12$, near the repelling fixed points $q_-$ and
$Jq_-$. The reflected images $h(U_\flat)$, $h=T^nI_0$, lie in the lower
half-plane.}
\label{fig:trapping}
\end{figure}
\begin{thm}[The full threshold]\label{thm:wall}
Let $36\le\tan^2A\le125/3$. Then $U_\flat$ satisfies \textup{(T1)}, \textup{(T2)} and \textup{(T3)}. Consequently $\rho_A$ is a discrete embedding.
\end{thm}

\begin{proof}
By Lemma~\ref{lem:strip}(3), $U_\flat\subset\mathcal W\subset\N$.

\textup{(T1)}. We have $I_2(U_\flat)=I_2(B\cap\mathcal W)\cup\mathcal O$. First, $U_\flat\cap\mathcal O\subset\mathcal W\cap\mathcal O=\emptyset$. Next, $B\cap I_2(B)=\emptyset$ by Lemma~\ref{lem:inversion}. Finally $\mathcal C\cap I_2(B\cap\mathcal W)=I_2(\mathcal O\cap B\cap\mathcal W)=\emptyset$.

\textup{(T2)}. The set $B\cap\mathcal W$ is open and contains $o$, by Lemma~\ref{lem:strip}(2).

\textup{(T3)} for $h=T^n$, $n\ne0$. By Lemma~\ref{lem:strip}(1), $h(U_\flat)\subset T^n(\mathcal W)\subset\mathcal O$. Hence $I_2h(U_\flat)\subset\mathcal C\subset U_\flat$.

\textup{(T3)} for $h=T^nI_0$, $n\in\Z$. We show that $I_2h(U_\flat)\subset B\cap\mathcal W$. The map $h$ acts on the horizontal coordinate by $z\mapsto-z-2n$, so it changes the sign of $\Im z$. Every point of $\cl B$ satisfies $|\Im z-\tfrac t2|\le r$, and
\[
0<r-\frac t2=\frac1{2\bigl(\sqrt{1+t^2}+t\bigr)}<\frac1{4t}\le\frac1{24}.
\]
Thus $\Im z\ge-\tfrac1{24}$ on $\cl B$ and $\Im z<\tfrac1{24}$ on $h(B)$. By Lemma~\ref{lem:caps}, $\Im z>\tfrac1{12}$ on $\mathcal C^{*}\supset\mathcal C$ and $\Im z<-\tfrac1{12}$ on $h(\mathcal C)$. Since $h(U_\flat)\subset h(B)\cup h(\mathcal C)$, we conclude
\[
h(U_\flat)\cap\mathcal C^{*}=\emptyset,\qquad h(\mathcal C)\cap\cl B=\emptyset .
\]
Also $h(B)\cap\cl B=\emptyset$, because $\rho(o,ho)\ge2r$ by Lemma~\ref{lem:centres} and the triangle inequality. Therefore $h(U_\flat)\cap\cl B=\emptyset$, and $I_2h(U_\flat)\subset B$ by Lemma~\ref{lem:inversion}. Finally, $h(U_\flat)\cap\mathcal C^{*}=\emptyset$ means $I_2h(U_\flat)\cap\bigl(\partial\CH\setminus\mathcal W\bigr)=\emptyset$, that is, $I_2h(U_\flat)\subset\mathcal W$.

Theorem~\ref{thm:D} now proves discreteness and faithfulness.
\end{proof}

\begin{proof}[Proof of Theorem~\ref{thm:GP}]
Proposition~\ref{prop:initial36} proves discreteness and faithfulness for $t^2\le36$. Theorem~\ref{thm:wall} proves them for $36\le t^2\le125/3$. Together they cover the full range, including the critical parameter. This proves the Goldman--Parker conjecture. 
\end{proof}

\appendix
\section{The exact Ford-ball threshold}\label{app:sharp-ball}

This appendix proves Lemma~\ref{lem:key} and hence the sharp statement in Theorem~\ref{thm:B}. It is independent of the proof of the full Goldman--Parker theorem, which uses only Proposition~\ref{prop:initial36} for the initial range.

\begin{proof}[Proof of Lemma~\ref{lem:key}]
Use the normalization $\Lambda$ and the distances $D_\pm$ from Proposition~\ref{prop:initial36}. Equation~\eqref{eq:ball-average} shows that any intersection of the closed balls forces $Q_t(y)\le0$. If one of the balls is open, it forces $Q_t(y)<0$. Conversely, at $X=v=0$ the two fourth-power distances are equal to
\[
(1+y^2)^2+4(t+y)^2.
\]
Thus the open balls intersect exactly when $\min_y Q_t(y)<0$; the closed balls intersect exactly when this minimum is nonpositive.

Now $Q_t'(y)=64(y^3+3y+2t)$. The cubic is strictly increasing, so the minimum is attained at the unique point $y=-Y$, where $Y\ge0$ satisfies $Y^3+3Y=2t$. At that point,
\[
t-Y=\frac12Y(1+Y^2),\qquad
1+t^2=\frac14(1+Y^2)^2(4+Y^2).
\]
Putting $\eta=Y^2$, we obtain
\[
\min Q_t=\frac{(1+\eta)^3}{16}
\bigl[256-(1+\eta)(4+\eta)^2\bigr].
\]
The function $(1+\eta)(4+\eta)^2$ is strictly increasing from $16$ to infinity for $\eta\ge0$, and $t=(Y^3+3Y)/2$ is strictly increasing in $Y$. Hence there is a unique threshold $t_*$, determined by $(1+\eta_*)(4+\eta_*)^2=256$. For $t<t_*$ the minimum is positive, for $t=t_*$ it is zero, and for $t>t_*$ it is negative.

At equality, equation~\eqref{eq:ball-average} forces $X=v=0$ and $y=-Y_*$. Thus the closed balls have exactly one common point. If $t>t_*$, this same choice $X=v=0$, $y=-Y$ lies inside both balls. If $t\le t_*$, no intersection is possible when one of the balls is open.

Finally put $w=\sqrt{1+\eta_*}$. The threshold equation becomes $w^2(w^2+3)^2=256$, or $w^3+3w=16$, with $w>0$. It follows that
\[
1+t_*^2=\frac14w^4(w^2+3)=4w^3=4(16-3w),
\]
so $t_*^2=63-12w$. Equivalently,
$w=\sqrt[3]{8+\sqrt{65}}-\sqrt[3]{\sqrt{65}-8}$.
\end{proof}

\medskip 

\subsection*{AI declaration}
AI-assisted tools were used in preparing the manuscript. The wall construction in Section~\ref{sec:wall} was proposed by Claude Fable 5.2 (Anthropic), which was also used to check the computations in the proof. All mathematical statements and proofs were independently reviewed and verified by the authors, who take full responsibility for the content of the paper.

\end{document}